\documentclass[]{amsart}

\usepackage{amsmath}

\theoremstyle{plain}
\newtheorem{thm}{Theorem}[section]
\newtheorem{prop}[thm]{Proposition}

\theoremstyle{definition}
\newtheorem{defn}[thm]{Definition}

\newtheorem{rmk}[thm]{Remark}

\newtheorem{quest}[thm]{Question}

\newcommand{\vol}{\ensuremath{\operatorname{vol}}}
\newcommand{\real}{\ensuremath{\operatorname{Re}}}
\newcommand{\imag}{\ensuremath{\operatorname{Im}}}

\newcommand{\G}{\ensuremath{\operatorname{G_2}}}
\newcommand{\SP}{\ensuremath{\operatorname{Spin (7)}}}
\newcommand{\SUn}{\ensuremath{\operatorname{SU} (n)}}
\newcommand{\SUt}{\ensuremath{\operatorname{SU} (2)}}
\newcommand{\SPo}{\ensuremath{\operatorname{Sp} (1)}}

\newcommand{\cf}{\ensuremath{\check{f}}}
\newcommand{\ph}{\ensuremath{\varphi}}

\newcommand{\Ph}{\ensuremath{\Phi}}

\newcommand{\wtwom}{\ensuremath{\wedge^2_{-}}}
\newcommand{\R}{\ensuremath{\mathbb R}}
\newcommand{\C}{\ensuremath{\mathbb C}}
\newcommand{\Q}{\ensuremath{\mathbb H}}
\newcommand{\PR}{\ensuremath{\mathbb P}}
\newcommand{\Oc}{\ensuremath{\mathbb O}}
\newcommand{\Qu}{\ensuremath{\mathbb H}}
\newcommand{\st}{\ensuremath{\ast}}
\newcommand{\hk}{\ensuremath{\lrcorner}}
\newcommand{\stph}{\ensuremath{\ast \varphi}}

\newcommand{\rest}[2]{\ensuremath{{#1} |_{{}_{#2}}}}
\newcommand{\nup}{\ensuremath{\nu^{\perp}}}
\newcommand{\Xs}[3]{\ensuremath{ {X({#1},{#2},{#3}) }^{\flat} }}
\newcommand{\ws}{\ensuremath{w^{\flat}}}
\newcommand{\us}{\ensuremath{u^{\flat}}}
\newcommand{\vs}{\ensuremath{v^{\flat}}}
\newcommand{\ys}{\ensuremath{y^{\flat}}}

\newcommand{\spm}{\ensuremath{{\, \slash \! \! \! \! \mathcal
S}_{\! \!{{}^-}}}}
\newcommand{\sppm}{\ensuremath{{\, \slash \! \! \! \! \mathcal
S}_{\! \!{{}^{\pm}}}}}
\newcommand{\spi}{\ensuremath{{\, \slash \! \! \! \! \mathcal S}}}

\numberwithin{equation}{section}
\numberwithin{table}{section}
\numberwithin{figure}{section}

\begin{document}

\title[Calibrated subbundles in non-compact manifolds]{Calibrated
subbundles in non-compact \\ manifolds of special holonomy}

\author{Spiro Karigiannis}
\thanks{Spiro Karigiannis was partially supported by NSERC}
\address{Mathematics Department\\ McMaster University}
\email{spiro@math.mcmaster.ca}
\author{Maung Min-Oo}
\thanks{Maung Min-Oo was partially supported by NSERC}
\address{Mathematics Department\\ McMaster University}
\email{minoo@mcmaster.ca}
\subjclass{53, 58} \date{\today}

\begin{abstract}
This paper is a continuation of~\cite{IKM}. We first construct special
Lagrangian submanifolds of the Ricci-flat Stenzel metric (of holonomy
\SUn) on the cotangent bundle of $S^n$ by looking at the conormal
bundle of appropriate submanifolds of $S^n$. We find that the
condition for the conormal bundle to be special Lagrangian is the same
as that discovered by Harvey-Lawson for submanifolds in $\R^n$ in
their pioneering paper~\cite{HL}. We also construct calibrated
submanifolds in complete metrics with special holonomy \G\ and \SP\
discovered by Bryant and Salamon~\cite{BS} on the total spaces of
appropriate bundles over self-dual Einstein four manifolds. The
submanifolds are constructed as certain subbundles over immersed
surfaces. We show that this construction requires the surface to be
minimal in the associative and Cayley cases, and to be (properly
oriented) real isotropic in the coassociative case. We also make
some remarks about using these constructions as a possible local
model for the intersection of compact calibrated submanifolds in a
compact manifold with special holonomy.
\end{abstract}

\maketitle

\section{Introduction} \label{introsec}

The study of calibrated geometries begun in the paper~\cite{HL} of
Harvey and Lawson. Calibrated submanifolds (in particular special
Lagrangian submanifolds) are believed to play a crucial role in mirror
symmetry~\cite{SYZ} and M-theory, and hence they have recently
received much attention. There has been extensive research done on
special Lagrangian submanifolds of $\mathbb C^n$, most notably by
Joyce but see also~\cite{J2} and the many references contained
therein. Much less progress has been made in studying associative,
coassociative, and Cayley submanifolds even in flat space. The
earliest explicit non-flat examples of special holonomy metrics were
constructed on vector bundles. These explicit metrics are all
cohomogeneity one examples and are obtained by reducing the conditions
for special holonomy to an exactly solvable ordinary differential
equation. Explicit Calabi-Yau metrics were found on the cotangent
bundle of sphere, initially discovered by Eguchi-Hanson for $S^2$ and
Candelas and others for $S^3$, but see Stenzel~\cite{St} for the
general case. Similarly Calabi discovered hyper-K\"{a}hler metrics on
the cotangent bundle of the complex projective space~\cite{Ca} and
Bryant and Salamon~\cite{BS} found explicit examples of metrics of
full holonomy \G\ and \SP\ on the bundles of anti-self-dual $2$-forms
and negative chirality spinors over specific four-manifolds. (See
Remark~\ref{orientationsrmk} for a note on orientation conventions.)
These bundles, although non-compact, also serve as local models for a
general metric of special holonomy and they have also received a lot
of attention from mathematical physicists, who have generalized these
metrics and studied them in detail~\cite{AW, CGLP1, CGLP2, GPP, GS}.

In our first paper~\cite{IKM}, along with Marianty Ionel we
generalized a bundle construction of Harvey and Lawson for special
Lagrangian submanifolds in $\mathbb C^n$ to analogous
constructions of coassociative, associative, and Cayley
submanifolds in $\mathbb R^7$ and $\mathbb R^8$. In this paper we
further generalize this construction to the case of several
explicit, non-flat, non-compact manifolds with complete metrics of
special holonomy which are vector bundles over a compact base. The
authors recommend that readers first consult~\cite{IKM}, as many
of the calculations, especially in Section~\ref{exceptionalsec}
are very similar and are covered in more detail in~\cite{IKM}. In
particular, without further mention, all of our local
calculations are done using normal coordinates.

In Section~\ref{calibreviewsec} we briefly review the relevant
facts from calibrated geometry that we will use, and set up some
notation. In particular, it should be noted that in
Propositions~\ref{assocprop} and~\ref{cayleyprop} we present
alternative characterizations of the associative and Cayley
conditions. These characterizations are entirely in terms of the
calibrating forms and the associated cross products and metrics
(which are all derivable from the forms). This is similar to the
special Lagrangian and coassociative conditions. In~\cite{IKM} our
proofs in the associative and Cayley cases relied on a choice of
identification of the tangent spaces with octonions or purely
imaginary octonions and was perhaps not as satisfying. At least the
invariant description of the Cayley condition seems not to have
appeared in the literature before.

In Section~\ref{stenzelsec} we describe the Stenzel Calabi-Yau
metrics on $T^* (S^n)$ and show that the conormal bundle over an
immersed submanifold $X$ in $S^n$ is special Lagrangian with
respect to some phase (which depends on the codimension of $X$ in
$S^n$) if and only if $X$ is {\em austere} in $S^n$. This is the
same result as Harvey and Lawson found~\cite{HL} for $\mathbb C^n$
but it is perhaps surprising, especially since the complex
structure on $T^* (S^n)$ is obtained in an extremely different way
from that of $\mathbb C^n = T^*(\mathbb R^n)$, namely by
identifying it with a complex quadric hypersurface in $\mathbb
C^{n+1}$.

In Section~\ref{exceptionalsec} we construct coassociative and
Cayley submanifolds in $\wedge^2_- (S^4)$ and $\wedge^2_- (\C
\PR^2)$ by taking vector subbundles over an immersed surface
$\Sigma$ in the base. As in~\cite{IKM}, the associative
construction requires $\Sigma$ to be minimal, while the
coassociative case needs $\Sigma$ to be (properly oriented)
isotropic. (Sometimes also called superminimal.) In this case it
is perhaps not so surprising that the results are the same as in
the flat case, since the calculations are extremely similar,
differing bascially only by the presence of some conformal scaling
factors. This is entirely due to the fact that these cohomogeneity
one metrics have a high degree of symmetry. We also construct
Cayley submanifolds in the negative spinor bundle $\spm (S^4)$
over $S^4$ by taking rank $2$ vector bundles over a minimal
surface $\Sigma$ in $S^4$. The result is again the same as the flat
case of $\mathbb R^8$ found in~\cite{IKM} although this time the
calculation is done in a very different way. It should also be
noted that in the case of $\mathbb R^8$, we obtained degenerate
examples. That is, they were products of lower order constructions.
However, this time in the case of $\spm(S^4)$ the Cayley examples
are not degenerate.

Finally in Section~\ref{intersectionsec} we make some remarks
about how these constructions might be used as local models for
the intersections of compact calibrated submanifolds of a compact
manifold with special holonomy. We hope to expand upon this topic
further in a subsequent paper.

{\em Remark.} Similar although different statements to some of the
results of Section~\ref{exceptionalsec} appeared, without proof, in an
unpublished preprint by S.H. Wang~\cite{W} back in 2001.  As remarked
in ~\cite{IKM}, the original statement which appeared in the preprint
was incorrect, but the authors were recently notified by Robert Bryant
that a corrected version of Wang's paper will appear soon.

{\em Acknowledgements.} The authors would like to thank Marianty
Ionel and Nai-Chung Conan Leung for helpful discussions.

\section{Review of Calibrated Geometries} \label{calibreviewsec}

In this section we review the necessary facts about the calibrated
geometries that we study in this paper, and set up notation. Some
references are~\cite{HL, J2,J3}. Calibrated submanifolds are a
distinguished class of submanifolds of a Riemannian manifold
$(M,g)$ which are absolutely volume minimizing in their homology
class. Being minimal is a second order differential condition, but
being calibrated is a {\em first order} differential condition.

\begin{defn} \label{calibdefn}
A closed $k$-form $\alpha$ on $M$ is called a calibration if it
satisfies $\alpha(e_1, \ldots, e_k) \leq 1$ for any choice of $k$
orthonormal tangent vectors $e_1, \ldots , e_k$ at any point $p
\in M$. A calibrated subspace of $T_p (M)$ is an oriented
$k$-dimensional subspace $V_p$ for which $\alpha (V_p) = 1$. Then
a calibrated submanifold $L$ of $M$ is a $k$-dimensional oriented
submanifold for which each tangent space is a calibrated subspace.
Equivalently, $L^k$ is calibrated if
\begin{equation*}
\rest{\alpha}{L} = \vol_L
\end{equation*}
where $\vol_L$ is the volume form of $L$ associated to the induced
Riemannian metric from $M$ and the choice of orientation.
\end{defn}

Here are the four main examples of calibrated geometries. (More will
be said below about \G\ and \SP\ structures.)

{\bf I.} Complex submanifolds $L^{2 k}$ (of complex dimension $k$) of
a {\em K\"ahler} manifold $M$ where the calibration is given by
$\alpha = \frac{\omega^k}{k!}$, and $\omega$ is the K\"ahler form
on $M$. K\"ahler manifolds are characterized by having Riemannian
holonomy contained in $\operatorname{U} (n)$, where $n$ is the
complex dimension of $M$. These submanifolds come in all even real
dimensions.

{\bf II.} Special Lagrangian submanifolds $L^n$ {\em with phase $e^{i
\theta}$} of a {\em Calabi-Yau} manifold $M$ where the calibration
is given by $\real (e^{i \theta} \Omega)$, where $\Omega$ is the
holomorphic $(n,0)$ volume form on $M$. Calabi-Yau manifolds have
Riemannian holonomy contained in $\SUn$. Special Lagrangian
submanifolds are always half-dimensional, but there is an $S^1$
family of these calibrations for each $M$, corresponding to the
$e^{i \theta}$ freedom of choosing $\Omega$. Note that Calabi-Yau
manifolds, being K\"ahler, also possess the K\"ahler calibration.

{\bf III.} Associative submanifolds $L^3$ and coassociative
submanifolds $L^4$ of a \G\ manifold $M^7$. Here the calibrations
are given by the $3$-form $\ph$ and the $4$-form $\stph$,
respectively, where $\ph$ is the fundamental $3$-form
corresponding to the \G -structure. \G\ manifolds have Riemannian
holonomy contained in \G. These calibrated submanifolds only come
in dimensions $3$ and $4$.

{\bf IV.} Cayley submanifolds $L^4$ of a \SP\ manifold $M^8$. Here the
calibration is given by the $4$-form $\Ph$ which is the fundamental
$4$-form corresponding to the \SP -structure. \SP\ manifolds have
Riemannian holonomy contained in \SP. These calibrated submanifolds
only come in dimension $4$.

\begin{rmk}
If $M^{4n}$ is a {\em hyper-K\"ahler manifold}, which means its
Riemannian holonomy is contained in $\operatorname{Sp}(n)$, then it
has an $S^2$ family of K\"ahler structures and each one is
Calabi-Yau. There is thus a wealth of calibrated submanifolds in the
hyper-K\"ahler case. Also, a Calabi-Yau manifold $M^8$ of complex
dimension $4$ is always a \SP\ manifold, and thus contains special
Lagrangian, complex, and Cayley submanifolds.
\end{rmk}

In practice, it is not easy to check if $\rest{\alpha}{L} = \vol_L$
but there are alternative, equivalent conditions for a submanifold to
be calibrated which we now describe.

{\bf I.} Complex submanifolds $L$ of a K\"ahler manifold $M$ are
characterized by the fact that their tangent spaces are invariant
under the action of the complex structure $J$ on $M$.

{\bf II.} Harvey and Lawson showed in~\cite{HL} that, up to a possible
change of orientation, $L$ is special Lagrangian of phase $e^{i
\theta}$ if and only if
\begin{eqnarray} \label{slagconditioneq1}
\rest{\omega}{L} & = & 0 \qquad \qquad \text{ and} \\
\label{slagconditioneq2} \rest{\imag( e^{i \theta} \Omega)}{L} & = & 0
\end{eqnarray}
Condition~\eqref{slagconditioneq1} say that $L$ is Lagrangian,
while~\eqref{slagconditioneq2} is the {\em special} condition.

{\bf III.} A Riemannian manifold $M^7$ which possesses a \G\ structure
has a globally defined, two-fold vector cross product
\begin{eqnarray*}
\times : & & T (M) \times T (M) \to T (M) \\
& & (v, w ) \mapsto v \times w
\end{eqnarray*}
which satisfies
\begin{eqnarray*}
v \times w = - w \times v & & \qquad \text{ $\times$ is alternating}
\\ \langle v \times w , v \rangle = 0 & & \qquad \forall \, v,w \text{
(orthogonal to its arguments)} \\ {| v \times w |}^2 = {|v \wedge
w|}^2 & & \qquad \forall \, v,w
\end{eqnarray*}
where $\langle \cdot, \cdot \rangle$ is the Riemannian metric on $M$
and $| \cdot |$ is its associated norm. The metric, cross
product, and fundamental $3$-form $\ph$ are related by
\begin{equation} \label{phimetriceq}
\ph (u, v, w) = \langle u \times v , w \rangle
\end{equation}
from which it follows that
\begin{equation} \label{cross7eq}
{(u \times v)}^{\flat} = v \hk u \hk \ph
\end{equation}
where ${}^{\flat}$ is the isomorphism from vector fields to one-forms
induced by the Riemannian metric. It is shown in~\cite{HL} that a
$3$-dimensional submanifold $L^3$ is associative if and only if
its tangent space is preserved by the cross product $\times$.
Similarly, a $4$-dimensional submanifold $L^4$ is coassociative
if and only if $u \times v$ is a normal vector for every pair of
vectors $u,v$ tangent to $L^4$. There exist vector valued
alternating $3$ and $4$-forms on $M$ called the associator and
coassociator which vanish on associative and coassociative
submanifolds, respectively, but these are difficult to work with
directly as they are related to octonion algebra. In~\cite{HL}
Harvey and Lawson showed that the coassociative condition is
equivalent (up to a change of orientation), to the vanishing of
the $3$-form $\ph$:
\begin{equation} \label{coassconditioneq}
\rest{\ph}{L^4} = 0
\end{equation}
This reformulation should be compared to~\eqref{slagconditioneq1}
and~\eqref{slagconditioneq2}.

We now present an alternative characterization of the associative
condition. Let $u,v,w$ be a linearly independent set of tangent
vectors at a point $p \in M$. We want to check when the
$3$-dimensional subspace that they span is an associative
subspace. Now if we have chosen an identification of $T_p M$ with
$\imag \Oc$, then we need to check the vanishing of the associator:
\begin{equation*}
[u, v, w] = u(vw) - (uv)w
\end{equation*}
When $u$ and $v$ are imaginary octonions, their product is $uv = -
\langle u, v \rangle + u \times v$, in terms of the inner product
and the cross product. Thus we have
\begin{eqnarray*}
[ u, v, w ] & = & u ( - \langle v, w \rangle + v \times w ) - ( -
\langle u, v \rangle + u \times v ) w \\ & = & - \langle v , w
\rangle u - \langle u, v \times w \rangle  + u \times (v
\times w) \\ & & {}+ \langle u, v \rangle w  + \langle u \times v
, w \rangle  - (u \times v ) \times w \\ & = & \langle u, v
\rangle w - \langle v, w \rangle u + u \times (v \times w) - (u
\times v) \times w
\end{eqnarray*}
where we have used~\eqref{phimetriceq} to cancel two of the terms.
Now from Lemma 2.4.3 in~\cite{K} we have the formula
\begin{equation*}
u \times (v \times w) = - \langle u, v \rangle w + \langle u, w
\rangle v - {(u \hk v \hk w \hk \stph)}^{\#}
\end{equation*}
Substituting this into the above expression for the associator and
simplifying, we obtain
\begin{equation*}
[ u, v, w ] = - 2 {(u \hk v \hk w \hk \stph)}^{\#} 
\end{equation*}
Thus we have proved the following:
\begin{prop} \label{assocprop}
The subspace spanned by the tangent vectors $u,v,w$ is an
associative subspace if and only if
\begin{equation} \label{assconditioneq}
u \hk v \hk w \hk \stph = 0.
\end{equation}
\end{prop}
\begin{rmk} \label{assconditionrmk}
The left hand side of~\eqref{assconditioneq} is (using the
metric isomorphism) a vector valued $3$-form which is invariant
under the action of \G. Therefore representation theory arguments
say it must be the associator, and here we show this directly. 
\end{rmk}

{\bf IV.} A Riemannian manifold $M^8$ which possesses a \SP\ structure
has a globally defined, three-fold vector cross product
\begin{eqnarray*}
X : & & T (M) \times T (M) \times T(M) \to T (M) \\
& & (u, v, w ) \mapsto X(u, v, w)
\end{eqnarray*}
which satisfies
\begin{eqnarray*}
X(u, v, w) & & \qquad \text{ is totally skew-symmetric} \\ \langle
X(u,v,w) , u \rangle = 0 & & \qquad \forall \, u,v,w \text{
(orthogonal to its arguments)} \\ {| X(u,v,w) |}^2 = {|u \wedge v
\wedge w|}^2 & & \qquad \forall \, u,v,w
\end{eqnarray*}
where $\langle \cdot, \cdot \rangle$ is the Riemannian metric on $M$
and $| \cdot |$ is its associated norm. As in the \G\ case, the
metric, cross product, and fundamental
$4$-form $\Ph$ are related by
\begin{equation} \label{Phimetriceq}
\Ph (u, v, w, y) = \langle X(u, v , w ), y \rangle
\end{equation}
from which it follows that
\begin{equation} \label{cross8eq}
{X(u, v, w)}^{\flat} = w \hk v \hk u \hk \Ph.
\end{equation}

It is shown in~\cite{HL} that a $4$-dimensional submanifold $L^4$ is
Cayley if and only if its tangent space is preserved by the cross
product $X$. As in the \G\ case, there exists a rank $7$ bundle
valued $4$-form $\eta$ on $M$ that vanishes on Cayley submanifolds.
This form $\eta$ is defined in terms of octonion multiplication.
Let $u,v,w,y$ be a linearly independent set of tangent vectors at a
point $p \in M$. We want to check when the $4$-dimensional
subspace that they span is Cayley subspace. Assuming an explicit
identification of $T_p M$ with $\Oc$, the form $\eta$ is:
\begin{equation*}
\eta = \frac{1}{4} \imag \left( \bar u X (v,w,y) + \bar v X (w,u,y)
+ \bar w X (u,v,y) + \bar y X (v,u,w) \right) 
\end{equation*}
We now describe a characterization of the Cayley
condition which is analogous to~\eqref{assconditioneq}, that does
not seem to have explicitly appeared in the literature before. The
fact we use is the following. The space of
$2$-forms on $M$ splits as $\wedge^2 = \wedge^2_7 \oplus
\wedge^2_{21}$, where at each point $\wedge^2_k$ is
$k$-dimensional. (see~\cite{J1,K}.) One can check by explicit
computation that if $u$ and $v$ are tangent vectors, identified as
octonions, then
\begin{equation*}
\imag(\bar u v) \cong \pi_7 (\us \wedge \vs)
\end{equation*}
where $\pi_7$ is projection onto $\wedge^2_7$. Thus, up to
isomorphism, the expression for the form $\eta$ becomes
\begin{equation*}
\eta = \pi_7 \left( \us \wedge \Xs{v}{w}{y} + \vs \wedge
\Xs{w}{u}{y} + \ws \wedge \Xs{u}{v}{y} + \ys \wedge \Xs{v}{u}{w}
\right)
\end{equation*}
We have an explicit formula for the projection $\pi_7$ in terms of
the $4$-form $\Ph$. (See~\cite{K}, for example, although we differ
by a sign here because of the opposite choice of orientation.)
This formula is
\begin{equation*}
\pi_7 ( \us \wedge \vs ) = \frac{1}{4} \left( \us \wedge \vs + u
\hk v \hk \Ph \right)
\end{equation*}
Combining these expressions, we have proved the following:
\begin{prop} \label{cayleyprop}
The subspace spanned by the tangent vectors $u,v,w,y$ is a 
Cayley subspace if and only if the $\wedge^2_7$ valued $2$-form
$\eta$ vanishes:
\begin{eqnarray*}
\eta = \us \wedge \Xs{v}{w}{y} + u \hk X(v,w,y) \hk \Ph + \vs
\wedge \Xs{w}{u}{y} + v \hk X(w,u,y) \hk \Ph & & \\ \nonumber {}+
\ws \wedge \Xs{u}{v}{y} + w \hk X(u,v,y) \hk \Ph + \ys \wedge
\Xs{v}{u}{w} + y \hk X(v,u,w) \hk \Ph = 0 & & 
\end{eqnarray*}
\end{prop}

\begin{rmk} \label{instantonsbranesrmk}
It should be evident that calibrated submanifolds seem to fall into
two different categories. There are those whose tangent spaces are
preserved by a cross product operation. These are the complex,
associative, and Cayley submanifolds, whose tangent spaces are
preserved by $J$, $\times$, and $X$, respectively. These are called
{\em instantons}. There are also those which are determined by the
vanishing of differential forms, namely the special Lagrangian and
coassociative submanifolds, and these are called {\em branes}. Branes
have a nice, unobstructed deformation theory, which was first
studied by McLean~\cite{M}. Instantons, on the other hand, are
generally obstructed and are more complicated to study.
See~\cite{LL} for more details on the differences between branes
and instantons.
\end{rmk}

\section{Special Lagrangians in $T^* (S^n)$ with the Stenzel metric}
\label{stenzelsec}

In this section we construct special Lagrangian submanifolds in
$T^* (S^n)$ with the Calabi-Yau metric discovered by
Stenzel~\cite{St} and discussed in detail in~\cite{CGLP1}.

It is a classical fact that if $X^p$ is a $p$-dimensional
submanifold of $\R ^n$, then the {\em conormal bundle} $N^* (X)$ is a
{\em Lagrangian} submanifold of the symplectic manifold $T^* (\R ^n)$,
with its canonical symplectic structure. Harvey and Lawson found
conditions~(\cite{HL}, Theorem III.3.11) on the immersion $X
\subset \R^n$ that makes $N^* (X)$ a special Lagrangian submanifold
of $T^* (\R^n) \cong \C^n$, in terms of the second fundamental
form of the immersion. We generalize this construction to the case
of the Calabi-Yau metric on $T^* (S^n)$, which we now describe.

Following Sz\"{o}ke~\cite{Sz}, we can map the space
\begin{equation*}
T^* (S^n) = \{ (x, \xi) \in \R^{n+1} \times \R^{n+1} | |x|=1, \langle
x,\xi \rangle = 0 \}
\end{equation*}
diffeomorphically and equivariantly with respect to
$\operatorname{SO}(n;\R) \subset \operatorname{O}(n;\C)$ onto the
complex quadric
\begin{equation*}
Q =  \{ (z_0, \ldots, z_n) \in \C^{n+1}  | \sum z_k^2 =1  \}
\end{equation*}
in $\C^{n+1}$ by
\begin{eqnarray} \nonumber
\Psi :  T^*S^n & \rightarrow &  Q \\ \label{szokemapeq}
(x, \xi) & \mapsto  &  x\, \cosh |\xi| + i \frac{\xi}{|\xi|} \sinh (|\xi|)
\end{eqnarray}
In this way $Q \cong T^* (S^n)$ inherits a complex structure, since
it is a complex hypersurface of $\C^{n+1}$. It also posseses a
holomorphic $(n,0)$ form $\Omega$ which is defined by
\begin{equation} \label{holovolformeq}
\Omega (v_1, \ldots, v_n) = \left( dz_0 \wedge dz_1 \wedge \ldots
\wedge dz_n \right) (Z, v_1, \ldots, v_n)
\end{equation}
where
\begin{equation*}
Z = z_0 \frac{\partial}{\partial z_0} + z_1
\frac{\partial}{\partial z_1} + \cdots + z_n
\frac{\partial}{\partial z_n}
\end{equation*}
is the holomorphic radial vector field on $\C^{n+1}$. With respect
to this complex structure, Stenzel showed~\cite{St} that there
exists a Ricci-flat K\"ahler metric on $T^* (S^n)$, thought of as
the quadic $Q$, whose K\"ahler form $\omega_{St}$, in a
neighbourhood of a point where $z_0 \neq 0$, is given by
\begin{equation} \label{stenzeleq1}
\omega_{St} = \frac{i}{2}\sum^{n}_{j,k = 1} a_{jk} dz_j \wedge d
\bar{z}_k
\end{equation}
where we have (see also Anciaux~\cite{An} for more details) that
\begin{equation} \label{stenzeleq2}
a_{jk} = \left( \delta_{jk} + \frac{ z_j \,\bar{z}_k}{|z_0|^2}
\right) u' + 2 \real \left( \bar{z}_j z_k -
\frac{\bar{z}_0}{z_0}z_jz_k \right) \,u''
\end{equation}
Here $u$ is a function of the radial variable $r = |z|$ and
satisfies a certain ordinary differential equation that makes the
metric Ricci-flat. The precise form of $u$ depends on the
dimension $n$ but it will not concern us since our results depend
only on the fact that $u$ is a function of $r$. We note that $r^2
= \cosh^2 |\xi| + \sinh ^2 |\xi|$. It is easy to check
from~\eqref{stenzeleq1} and~\eqref{stenzeleq2} that when
restricted to the zero section, this gives the standard round
metric on $S^n$. In dimension
$n=2$ this metric coincides with the well known Eguchi-Hanson and
Calabi metrics on $T^* (S^2)$.~\cite{Ca, CGLP1, CGLP2, DS}

Now let $X$ be a $p$-dimensional submanifold of the standard round
sphere $S^n$ with the induced metric. The conormal bundle of $X^p$
in $S^n$ will be denoted by $L = N^* (X) \subset T^*(S^n)$. Then
$L$ is a submanifold of dimension $n$ and can be locally be
parametrized as:
\begin{equation*}
(s,t) \mapsto \left( x(s), \Sigma t_k\nu^k \right) \quad \quad
s=(s_1, \ldots, s_p) ,\quad t=(t_{p+1}, \ldots, t_n)
\end{equation*}
where $x=(x_0, \ldots, x_n) \in X \subset S^n $ and
$\nu=(\nu^{p+1}, \ldots, \nu^n) \in \R^{n+1}$ are orthonormal
conormal vectors in $N^* (X)$. Let $e_1, \ldots ,e_p$ be an
orthonormal base of tangent vectors to $X$. Then $(e_0 = x(s),
e_1, \ldots, e_p,\nu^{p+1}, \ldots, \nu^n)$ form an adapted
orthonormal moving frame of $\R^{n+1}$ along the submanifold $X$.

We restrict the map in~\eqref{szokemapeq} to the subbundle $L = N^*
(X)$:
\begin{equation*}
\Psi \left( x(s), \Sigma t_k\nu^k \right) = x(s) \cosh |t| + i
\hat{\nu}(s,t) \sinh |t|
\end{equation*}
where $|t|^2 = t^2 _{p+1} + \ldots + t_n^2$, and $\hat{\nu}=
\frac{\Sigma t_k \nu^k}{|t|} $ is a unit conormal vector. Note that
$\hat{\nu}$ is homogeneous as a function of $t$. That is,
$\hat{\nu}(s, \lambda t) = \hat{\nu}(s,t) $ for all $\lambda \neq 0$
and $ \hat{\nu}(s,t) \sinh |t| $ is well defined for $t=0$.

\begin{thm} \label{slagthm}
The conormal bundle $L$ of a submanifold $X \subset S^n$ is special
Lagrangian in $T^* (S^n)$ equipped with the Ricci-flat Stenzel
metric if and only if $X$ is {\em austere} in $S^n$.
\end{thm}
\begin{proof}
We show that the tangent space of $L$ at each point is a
special Lagrangian subspace. Fix a point $(x, \xi) \in
L$. By the equivariance of the embedding we can choose an
orthonormal basis
$(e_0, \ldots, e_n)$ of $\R^{n+1}$ so that at the point $(x,
\xi)$ the moving frame is given by these vectors and so the
point has coordinates $(x(0) = e_0, \Sigma t_k \nu^k)$ with
$\nu^k(0) = e_k$, for $k = p+1, \ldots, n$. In fact, since we
still have the freedom of rotating the conormal vectors, we can
assume that $\hat{\nu} = \nu^{p+1} = e_{p+1}$. In other words, we
can rotate so that the point we are considering has $t$
coordinates $t_{p+1} = |t| = t \geq 0$ and $t_k = 0$ for $k = p+2,
\ldots, n$.

Now we compute a basis for the tangent space at this point
$\Psi(x,\xi) = e_0 \, \cosh|t| + i e_{p+1} \, \sinh|t| $. We
differentiate the immersion with respect to the $s$ and $t$
coordinates and evaluate at the point. From $s_1, \ldots, s_p$ we
have
\begin{equation} \label{stenzelbasis1eq}
E_j = \cosh |t|\, e_j + i \sinh |t| \, A^{\hat \nu} ( e_j ) \qquad
j = 1, \ldots, p
\end{equation}
where $A^{\hat{\nu}}$ is the second fundamental form in the
direction of the unit normal vector $\hat{\nu}$ of the submanifold
$X$ in $S^n$. That is, $A^{\hat{\nu}}(u) = \overline{\nabla}_u
\hat{\nu}$, where $\overline{\nabla}$ is the Levi-Civita
connection for the standard round metric on $S^n$. When we
differentiate with respect to $t_k$ we get
\begin{equation*}
F_k = x(s) \frac{\sinh |t|}{|t|} t_k + i \left( \nu^k
\frac{\sinh |t|}{|t|} + \left(\sum_l t_l \nu^l \right) \left(
\frac{|t|\cosh|t| - \sinh|t|}{|t|^3}t_k \right) \right)
\end{equation*}
Now we evaluate at our fixed point by putting $s=0$, $t_k=0$ for $k
\neq p+1$, and $t_{p+1} = |t|$ to obtain
\begin{eqnarray} \label{stenzelbasis2eq}
F_{p+1} & = & \sinh|t| \, e_0 + i \, \cosh|t| \, e_{p+1} \\
\nonumber F_k & = & \qquad \qquad i \, \frac{\sinh|t|}{|t|} \, e_k
\qquad k = p+2, \ldots , n
\end{eqnarray}

At the point $e_0 \, \cosh|t| + i e_{p+1} \, \sinh|t| $, $z_0 = \cosh
|t| \neq 0, \; z_{p+1} = i \, \sinh|t|$ and all the other coordinates
$z_1, \ldots, z_p, z_{p+2}, \ldots, z_n$ are zero. This simplifies (and
in fact diagonalizes) the Stenzel metric in~\eqref{stenzeleq1}
and~\eqref{stenzeleq2} and we have at that point
\begin{eqnarray*}
a_{jk} & = &  u'   \qquad j,k  \neq  p+1 \\ 
a_{p+1, p+1} & = & \big(1+ \tanh^2|t| \big) u'+ 4 \sinh^2|t| \, u''  
\end{eqnarray*}
and so 
\begin{equation*}
\omega_{St} = u' \, \frac{i}{2}\sum_{k=1}^{n} dz_k \wedge d \bar{z}_k +
\frac{i}{2} \big(u' \tanh^2|t| \, u'+ 4 u'' \sinh^2|t| \big) dz_{p+1}
\wedge d \bar{z}_{p+1}
\end{equation*}

Since from~\eqref{stenzelbasis1eq} the $E_j$'s have a zero
component in the $e_{p+1}$-direction,
$dz_{p+1} \wedge d \bar{z}_{p+1}$ vanishes on $E_j \wedge E_k$ for all
$j,k$ and we have
\begin{eqnarray*}
\omega_{St}(E_j,E_k) &=& u' \, \sinh|t| \cosh|t| \left( \langle
A^{\hat{\nu}}(e_j, e_k) \rangle - \langle A^{\hat{\nu}}(e_k, e_j)
\rangle \right) \\ & = & 0
\end{eqnarray*}
since the second fundamental form is symmetric.
From~\eqref{stenzelbasis1eq} and~\eqref{stenzelbasis2eq} we
see that $E_j$ has non-zero components only in the $z_1, \ldots
,z_p$ directions and
$F_k$ for $k= p+2, \ldots , n$ has a non-zero component only in
the $z_k$ direction. Hence
\begin{equation*}
\omega_{St}(E_j, F_k) = 0 \qquad \qquad j=1, \ldots , p \quad
\text{ and } \quad k=p+2, \ldots , n
\end{equation*}
Similarly, $F_{p+1}$ has non-zero components only in the direction
of $z_0$ and $z_{p+1}$. Thus
\begin{eqnarray*}
\omega_{St}(E_j,F_{p+1}) & = & 0  \\
\omega_{St}(F_k,F_{p+1}) & = & 0
\end{eqnarray*}
Thus we have shown that that $L = N^* (X)$ is always Lagrangian
with respect to the symplectic form associated to the Stenzel
metric for any submanifold $X$ of $S^n$.

In order to find the conditions for $L$ to be special Lagrangian, we
have to evaluate the holomorphic $(n,0)$-form $\Omega$ on the
tangent vectors $E_j$ and $F_k$ of our submanifold. In a
neighbourhood of a point where $z_0 \neq 0$, it follows
from~\eqref{holovolformeq} that
\begin{equation} \label{holovolformeq2}
\Omega = \frac{1}{2z_0}\, dz_1 \wedge \cdots \wedge dz_n
\end{equation}
This calculation is very similar to the original calculation done
by Harvey and Lawson~\cite{HL}, except that we have factors
involving the function $u$ and the hyperbolic trigonometric
functions of the radial variable $|t|$.

We can choose $e_1, \ldots, e_p$ to diagonalize the second fundamental
in the direction $\hat{\nu}$ at the point under consideration. Let
$\lambda_j$ be the corresponding eigenvalues (principal
curvatures). Then we have
\begin{eqnarray*}
E_j & = & \cosh|t| \, e_j + i \, \lambda_j \sinh|t| \, e_j \qquad j =
1, \ldots , p \\ F_{p+1} & = & \sinh|t| \, e_0 + i \, \cosh|t| \,
e_{p+1} \\ F_k & = & \qquad \qquad i \, \frac{\sinh|t|}{|t|} \, e_k
\qquad k = p+2, \ldots , n
\end{eqnarray*}
and hence, plugging into~\eqref{holovolformeq2}, 
\begin{eqnarray*}
& & \Omega(E_1 \wedge \cdots \wedge E_p \wedge F_{p+1} \wedge
\cdots \wedge F_{n}) \\ & = & \frac{1}{2 \cosh|t|} \cosh|t| {\left(
\frac{\sinh|t|}{|t|} \right)}^{n-p-1} i^{n-p} \prod_{j=1}^{p}
\left( \cosh|t| + i \lambda_j \sinh|t| \right) \\ & = & (***)
i^{n-p}  \prod_{j=1}^{p} \left( 1 + i \lambda_j \tanh|t|
\right)
\end{eqnarray*}
where $(***)$ denotes an always positive factor.
Hence from~\eqref{slagconditioneq2} we see that $L$ will be special
Lagrangian with phase $i^{p-n}$ if the product on the right hand
side above vanishes for all $t$. This happens if and only if all
odd symmetric polynomials in the eigenvalues $\lambda_j$ have to
be zero, or equivalently if all eigenvalues occur in pairs of
opposite signs.  This has to be true in all normal directons $\nu$
and so the submanifold must be {\em austere} as defined by Harvey
and Lawson~\cite{HL}. This completes the proof.
\end{proof}

\begin{rmk} \label{austerermk}
The first symmetric polynomial is the trace, so the submanifold
$M^p$ is necessarily minimal. If $p = 1,2$ this is the only
condition, but for $p \geq 3$ the austere condition is much
stronger than minimal.
\end{rmk}
\begin{rmk} \label{slagintersectionsrmk}
It is interesting to note that we cannot construct special
Lagrangian submanifolds in this way of arbitrary phase. The factor
of $i^{p-n}$ means that the allowed phase (up to orientation)
depends on the codimension $n-p$ of the immersion. We will say
more about this in Section~\ref{intersectionsec}.
\end{rmk}

Austere submanifolds have been studied for example in~\cite{Br3,
DF}. A particularly simple (and in some sense trivial) example
comes from equators: a sphere $S^p$ immersed in $S^n$ as an equator
is totally geodesic, and hence the conormal bundle $N^* (S^p)$ is
a special Lagrangian submanifold of $T^* (S^n)$ with respect to the
Stenzel metric. (Of phase $i^{n-p}$.)

\section{Calibrated submanifolds for the Bryant-Salamon metrics}
\label{exceptionalsec}

In this section we will construct calibrated submanifolds as
subbundles inside the Bryant-Salamon metrics~\cite{BS} of exceptional
holonomy \G\ or \SP\ which are themselves defined on appropriate
bundles over four manifolds with a self-dual Einstein metric.  The
subbundles are defined exactly in the same way as in [IKM], except
that the ambient manifold, instead of being flat $\R^7$ or $\R^8$ is
the total space of a vector bundle over a four manifold $X^4$.

\subsection{Calibrated submanifolds of $\wedge^2_- (X^4)$}
\label{g2sec}

Let $(X^4,g)$ be an oriented self-dual Einstein manifold. The
examples for which Bryant and Salamon obtained complete \G\
metrics are those with positive scalar curvature: $\C\PR^2$ and
$S^4$. Let $M^7 = \wtwom(T^*X^4)$ be the bundle of anti-self-dual
$2$-forms on $X^4$. This vector bundle has a connection induced by
the Levi-Civita connection of $(X,g)$. The tangent space
$T_{\omega}M$ of $M$ at a point $\omega \in \wtwom $ has therefore a canonical
splitting $T_{\omega}M \cong \mathcal{H}_{\omega }\oplus
\mathcal{V}_{\omega }$ into horizontal and vertical subspaces.

The projection map is a submersion and maps the horizontal space
isometrically onto the tangent space of the base manifold at that
point. The metric $g$ on the base $X^4$ has a unique lift to the
horizontal space $g_{\mathcal{H}}$. The vertical space
$\mathcal{V}_{\omega}$, which can be identified with the vector
space (the fibre) $\wtwom(T^*_{x}X)$ also has a natural metric
$g_{\mathcal{V}}$ induced by $g$.
 
\begin{thm} (Bryant-Salamon~\cite{BS})
There exist positive functions $u$ and $v$, depending only on the
radial coordinate in the vertical fibres and satisfying a certain set
of ordinary differential equations such that the metric
\begin{equation} \label{BS7metriceq}
g_{M^7} = u^2\,g_{\mathcal{H}} \oplus v^2 \, g_{\mathcal{V}}
\end{equation}
on the total space $M^7 = \wtwom(T^*X^4) $ of a
self-dual Einstein 4-manifold has $G_2$-holonomy with fundamental
$3$-form $\ph$ given by
\begin{equation*}
\ph = v^3 \vol_{\mathcal{V}} + u^2 v \, d\theta
\end{equation*}
where $\theta$ is the canonical (soldering) $2$-form on
$\wtwom(T^*X^4)$ and $\vol_{\mathcal{V}}$ is the volume $3$-form of
$g_{\mathcal{V}}$ on the vertical fibres.
\end{thm}
\begin{rmk}
The canonical $p$-form $\theta$ on $\wedge^p(T^*X)$ for any manifold
$X$ is defined to be $\theta(u_1\wedge \cdots \wedge u_p)_{\omega} =
\omega(\pi_*u_1\wedge \cdots \wedge \pi_*u_p)$, at the point $\omega$
where $\pi$ is the projection onto the base manifold.  For $p=1$ this
is the usual canonical $1$-form on $T^* (X)$.
\end{rmk}

Let $e^0, e^1, e^2, e^3$ be an orthonormal coframe for $T^* (X)$
and $f^1, f^2, f^3$ be an (orthonormal) basis of anti-self-dual
$2$-forms in the vertical fibres defined by $f^i = e^0 \wedge e^i
- e^j \wedge e^k$ with $i,j,k$ forming a cyclic permutation of
$1,2,3$. We denote horizontal lifts of tangent vectors $e_i$ on
the base to $\mathcal{H}$ by $\bar{e}_i$, with dual horizontal
$1$-forms $\bar{e}^i$. Similarly we think of the anti-self dual
two forms $f^i$ as being vertical tangent vectors $\cf^i$ in
$\mathcal{V}$ on the total space with dual vertical $1$-forms
$\cf_i$. Then locally the fundamental three form
$\ph$ is given by
\begin{eqnarray} \label{BS7phieq}
& \ph & = v^3 \left( \cf_1 \wedge \cf_2 \wedge \cf_3 \right) + u^2 v \,
\cf_1 \wedge ( \bar{e}^0 \wedge \bar{e}^1 - \bar{e}^2 \wedge
\bar{e}^3) \\ \nonumber & & {}+ u^2 v \, \cf_2 \wedge ( \bar{e}^0 \wedge
\bar{e}^2 - \bar{e}^3 \wedge \bar{e}^1) + u^2 v \, \cf_3 \wedge (
\bar{e}^0 \wedge \bar{e}^3 - \bar{e}^1 \wedge \bar{e}^2)
\end{eqnarray}
In this basis, the dual $4$-form is given by
\begin{eqnarray} \label{BS7stphieq}
& \stph & = u^4 \left( \bar{e}^0 \wedge \bar{e}^1 \wedge \bar{e}^2
\wedge \bar{e}^3 \right) - u^2 v^2 \, \cf_2 \wedge \cf_3 \wedge (
\bar{e}^0 \wedge \bar{e}^1 - \bar{e}^2 \wedge \bar{e}^3) \\ \nonumber
& & {}- u^2 v^2 \, \cf_3 \wedge \cf_1 \wedge ( \bar{e}^0 \wedge \bar{e}^2
- \bar{e}^3 \wedge \bar{e}^1) - u^2 v^2 \, \cf_1 \wedge \cf_2 \wedge (
\bar{e}^0 \wedge \bar{e}^3 - \bar{e}^1 \wedge \bar{e}^2)
\end{eqnarray}

It was proved in~\cite{BS} that the functions $u$ and $v$ are
globally defined and the Bryant-Salamon metric is complete only in
the cases where $X$ is either $S^4$ or $\C\PR^2$ with the standard
metrics (round metric on $S^n$ and Fubini-Study metric on
$\C\PR^2$.) In other cases, like for hyperbolic space, the functions
are not globally defined and we only obtain an incomplete metric
defined near the zero section of the vector bundle
$\wtwom(T^*X^4)$. Our constructions below of associative and
coassociative submanifolds are of a general nature and hence works
in both cases (complete or incomplete).

An oriented surface $\Sigma^2 \subset X^4$ equipped with the induced
metric defines a canonical lift
\begin{equation*} 
f^1_{\Sigma}: \Sigma^2 \longrightarrow M^7 = \wtwom(X^4)
\end{equation*}
locally defined by the anti-self-dual $2$-form $f^1 = e^1 \wedge
e^2 - \nu^1 \wedge \nu^2$, where $e^1, e^2$ are orthonormal co-tangent
vectors and $\nu^1, \nu^2$ are orthonormal conormal vectors to the
surface $\Sigma$. That is, $(e^1, e^2, \nu^1, \nu^2)$ is an
oriented adapted co-frame along the surface. It is easily seen
that $f^1_{\Sigma}$ is globally well defined and is independent of
the local frame. More invariantly we can define it by
\begin{equation*}
f^1_{\Sigma} = \vol_{\Sigma} - \st \vol_{\Sigma}
\end{equation*}
where $\vol_{\Sigma}$ is the induced volume form on $\Sigma$ and $\st$
is the Hodge star operator on $X^4$. The span of $f^1$ defines a
line bundle $L^3 \subset M^7 = \wtwom(X)$. We also define $L^{\perp}=
\{ \omega \in \wtwom \,| \, \omega \perp \omega^1 \}$ to be the (real)
two-dimensional subbundle orthogonal to $L$ with respect to the
Bryant-Salamon metric. Locally $L^{\perp}$ is spanned by the two
anti-self-dual 2-forms
\begin{equation*}
f^2= e^1 \wedge \nu^1 - \nu^2 \wedge e^2 \qquad \qquad f^3= e^1
\wedge \nu^2 - e^2 \wedge \nu^1
\end{equation*}

We want to determine necessary and sufficient conditions on the second
fundamental form of $\Sigma$ for $L$ to be associative and $L^{\perp}$
to be coassociative with respect to the Bryant-Salamon
$G_2$-structure on $M^7$.

\begin{thm}
The bundle $L$ defined above which is canonically associated to a
surface $\Sigma$ in a four dimensional self-dual Einstein manifold
$(X^4,g)$ is associative in $M^7 = \wtwom(T^*X)$ equipped with
the $\G$ metric of Bryant and Salamon if and only if $\Sigma$ is a
minimal surface in $X^4$. The bundle $L^{\perp}$ is coassociative
if and only if $\Sigma$ is a (propertly oriented) real
isotropic surface in $X^4$.
\end{thm}
\begin{proof}
We check that at each point, the tangent space is a calibrated
subspace. We begin with the associative case. At a point $t_1 f^1 \in L
$, the following three vectors form a basis of the tangent space $T_{t_1
f^1}$ of $L$. (We denote the dual vectors with a lower index.)
\begin{eqnarray}
\label {assocbasiseq} E_i & = & \bar{e_i} + t_1 \alpha(e_i, f^1)
\qquad i = 1, 2 \\ \nonumber F_1 & = & \cf^1
\end{eqnarray}
where the bar denotes the horizontal lift and $\alpha(e_i, f^1) =
\left( \overline{\nabla}_{e_{i}} f^1 \right)_{\mathcal{V}}$ is a
vertical vector which can be expressed (locally) in terms of the
second fundamental form of the submanifold as follows:
\begin{equation*}
\alpha(e_i, f^1) = \left( - A^{\nu_1}(e_i,e_1) -
A^{\nu_2}(e_i,e_2) \right) \cf_3 + \left( - A^{\nu_1}(e_i,e_2) +
A^{\nu_2}(e_i,e_1) \right) \cf_2
\end{equation*}
where we use the notation: $A^{\nu}(u,v) = \langle
\overline{\nabla}_{u} \nu, v \rangle = -\langle
\overline{\nabla}_{u}v, \nu \rangle$ for $u,v \in T (X)$ and $\nu \perp
T (X)$.

From Proposition~\ref{assocprop} we have that the subbundle $L$ is
associative if and only if the $1$-form $E_1 \hk E_2 \hk F_1 \hk
\stph$ vanishes at all points of $L$. Using~\eqref{assocbasiseq} and
~\eqref{BS7stphieq} we compute:
\begin{equation*}
F_1 \hk \stph = - u^2 v \left( \cf_2 \wedge ( \bar{e}^1 \wedge
\bar{\nu}^2 - \bar{e}^2 \wedge \bar{\nu}^1) - \cf_3 \wedge
( \bar{e}^1 \wedge \bar{\nu}^1 - \bar{\nu}^2 \wedge \bar{e}^2)
\right)
\end{equation*}
where $u,v$ are just functions. Using the symmetry of the second
fundamental form $A$ and the index notation $A^i_{jk} =
A^{\nu_i}(e_j,e_k)$, we continue to compute:
\begin{eqnarray*}
E_2 \hk F_1 \hk \stph & = & -u^2 v \left(  \cf_2 \wedge \bar{\nu}^1 +
\cf_3 \wedge \bar{\nu}^2 + t_1 (A^1_{12} + A^2_{22})(\bar{e}^1 \wedge
\bar{\nu}^1 - \bar{\nu}^2 \wedge \bar{e}^2) \right) \\ & & {}- u^2 v
\left( t_1 (-A^1_{22} + A^2_{12})(\bar{e}^1 \wedge \bar{\nu}^2 -
\bar{e}^2 \wedge \bar{\nu}^1) \right)
\end{eqnarray*}
and further
\begin{eqnarray*}
E_1 \hk E_2 \hk F_1 \hk \stph & = & -u^2 v \left( t_1 (A^1_{12} +
A^2_{22}) \bar{\nu}^1 + t_1 (-A^1_{22} + A^2_{12}) \bar{\nu}^2 \right)
\\ & & {} -u^2 v \left( t_1 (-A^1_{11} - A^2_{12}) \bar{\nu}^2 + t_1
(-A^1_{12} + A^2_{11}) \bar{\nu}^1 \right) \\ & = & -t_1 u^2 v \left(
(A^2_{11} + A^2_{22}) \bar{\nu}^1 - ( A^1_{11} + A^1_{22} )
\bar{\nu}^2 \right)
\end{eqnarray*}
Since $u$, $v$, are positive functions and since this expression must
vanish at all points on $L$ (that is, for all $t_1$), we must have
$A^1_{11} + A^1_{22} = 0$ and $A^2_{11} + A^2_{22} = 0$. Thus $L$
is associative if and only if $\Sigma$ is a minimal surface in
$X^4$, proving the first half of the theorem.

We now move on to the coassociative case. For the subbundle
$L^{\perp}$ we have the following description of a basis of four
tangent vectors at a given point $f = t_2 f^2 + t_3 f^3 $:
\begin{eqnarray*}
E_i & = & \bar{e_i} + t_2 \alpha(e_i, f^2) + t_3 \alpha(e_i,
f^3) \qquad i = 1, 2 \\ F_j & = &  \cf^j \qquad \qquad
\qquad \qquad \qquad \qquad \quad j = 2, 3
\end{eqnarray*}
Here the vertical correction terms are given by:
\begin{eqnarray*}
\alpha(e_i, f^2) & = & \left( \bar{\nabla}_{e_{i}} f^2
\right)_{\mathcal{V}} = \left( A^{\nu_1}(e_i,e_2) -
A^{\nu_2}(e_i,e_1) \right) \cf^1 \\ \alpha(e_i, f^3) & = & \left(
\bar{\nabla}_{e_{i}} f^3 \right)_{\mathcal{V}} = \left( 
A^{\nu_2}(e_i,e_2) + A^{\nu_1}(e_i,e_1) \right) \cf^1
\end{eqnarray*}

In order to check coassociativity, by~\eqref{coassconditioneq} we
need to check that $\rest{\ph}{L^{\perp}} = 0$. As in~\cite{IKM} we
define $\nu = t_2 \nu_1 + t_3 \nu_2$ and $\nup = - t_3 \nu_1 + t_2
\nu_2$ and thus
\begin{eqnarray*}
E_1 & = & \bar{e}_1 + \left( A^{\nu}_{12} - A^{\nup}_{11} \right)
\cf^1 \\ E_2 & = & \bar{e}_2 + \left( A^{\nu}_{22} -
A^{\nup}_{12} \right) \cf^1
\end{eqnarray*}
It is easy to compute that
\begin{eqnarray*}
\ph (E_1, E_2, \cdot) & = & E_2 \hk E_1 \hk \ph \\ & = & u^2 v
\left( \cf_1 + \left( \cdots \right) \bar{e}^1 + \left( \cdots
\right) \bar{e}^2 \right)
\end{eqnarray*}
and hence since $F_j = \cf^j$ we see that $\ph(E_1, E_2, F_2) =
\ph(E_1, E_2, F_3) = 0$ always. It remains to check when $\ph
(F_2, F_3, E_j) = 0$ for $j = 1,2$. Since $\ph(F_2, F_3, \cdot) =
v^3 \cf_1$ and $v$ is always positive, these become the conditions
\begin{equation} \label{coassresults1}
A^{\nu}_{12} - A^{\nup}_{11} = 0 \qquad \qquad 
A^{\nu}_{22} - A^{\nup}_{12} = 0
\end{equation}
for the tangent space at $(\mathbf x_0, t_2, t_3)$ to be
coassociative. We get two more conditions that must be satisfied by
demanding that the tangent space at $(\mathbf x_0, -t_3, t_2)$
also be coassociative. This corresponds to changing $t_2 \mapsto
-t_3$ and $t_3 \mapsto t_2$ in the above equations, which is
equivalent to $\nu \mapsto \nup$ and $\nup \mapsto -\nu$. This
gives
\begin{equation} \label{coassresults2}
A^{\nup}_{12} + A^{\nu}_{11} = 0 \qquad \qquad 
A^{\nup}_{22} + A^{\nu}_{12} = 0
\end{equation}
Conditions~\eqref{coassresults1} and~\eqref{coassresults2} are
exactly the same as those obtained in the case of $\mathbb R^7$
in~\cite{IKM}. These surfaces are called isotropic (with negative
orientation) or superminimal surfaces. These surfaces are
necessarily minimal, but the condition is in fact stronger (and
overdetermined). See~\cite{Br3, ES, IKM, Sa2} and the references
contained therein for more details.
\end{proof}
\begin{rmk}
Although the associative case is computed using a different method
from that of~\cite{IKM}, the calculations here and in
Section~\ref{cayleysec} are very similar to~\cite{IKM}, basically
differing by the presence of certain conformal scaling factors.
This is due to the high degree of symmetry in the cohomogeneity
one metrics.
\end{rmk}

\subsection{Cayley Submanifolds of $\spm(S^4)$}
\label{cayleysec}

In order to construct Cayley submanifolds, we now look at the
Bryant-Salamon construction on the negative spin bundle of four
manifolds.  Let $(X^4,g)$ be an oriented self-dual Einstein {\em spin}
manifold of positive scalar curvature. The only example now is $S^4$,
since $\C\PR^2$ is not spin. Let $M^8= \spm(X^4) \longrightarrow X$ be
the complex two-dimensional vector bundle of negative chirality
spinors on $S^4$.  This is in fact the quaternionic Hopf bundle of the
quaternionic projective line $\Qu \PR^1 \cong S^4$. Its unit sphere
bundle $S^7 \longrightarrow S^4$, can be viewed as the associated
principal $\SPo \cong \SUt$-bundle.  Note that $\operatorname{Spin}(4)
\cong \SUt \times \SUt$. This vector bundle has a natural Hermitian
inner product and a connection induced by the Levi-Civita connection
of the standard metric on $S^4$. The tangent space $T_{s}M$ of $M$ at
a point $ s \in \spm$ has therefore a canonical splitting $T_{s}M
\cong \mathcal{H}_{s }\oplus \mathcal{V}_{s }$ into horizontal and
vertical subspaces. It is well known that this connection defines the
standard $\SUt$-instanton on $S^4$ with (anti-) self-dual
curvature. The horizontal space of the connection is orthogonal to the
vertical space with respect to the standard metric on $S^7$ and the
curvature, which is the Lie bracket of horizontal vector fields
identifies the anti-self-dual $2$-forms on the base with the vertical
fibres which form the Lie algebra $su(2) \cong \R^3$. The projection
map is a submersion and maps the horizontal space isometrically onto
$T (S^4)$. The vertical space $\mathcal{V}_{s}$ also has a natural
induced metric $g_{\mathcal{V}}$ and the connection form is an
isomorphism between anti-self-dual $2$-forms and the Lie algebra of
$\SUt$.
 
\begin{thm} (Bryant-Salamon~\cite{BS})
There exist positive functions $u$ and $v$, depending only on the
radial coordinate in the vertical fibres and satisfying a certain
set of ordinary differential equations such that the metric
\begin{equation} \label{BS8metriceq}
g_{M^8} = u^2\,g_{\mathcal{H}} \oplus v^2 \, g_{\mathcal{V}}
\end{equation}
on the total space $M^8= \spm(S^4)$ has $\SP$-holonomy with
self-dual fundamental $4$-form $\Ph$ given by
\begin{equation*}
\Ph = u^4 \vol_{\mathcal{H}} + \,u^2v^2 \, \beta + v^4
\vol_{\mathcal{V}}
\end{equation*}
where $\vol_{\mathcal{H}}$, $\,\vol_{\mathcal{V}}$ are the volume
4-forms of $g_{\mathcal{V}}$, $\, g_{\mathcal{V}}$ on the
horizontal and vertical spaces respectively and $\beta$ is the
4-form defined as follows:
\begin{equation*}
\beta = \sum_{k=1}^{3} \omega_{k} \wedge \sigma^{k}
\end{equation*}
where $ \omega_{k}$ is an orthonormal basis for anti-self-dual
$2$-forms on the horizontal space and $ \sigma^{k}$ is the
corresponding orthonormal basis for anti-self-dual $2$-forms on the
vertical space.
\end{thm}
\begin{rmk}
Given an orthonormal basis of three anti-self-dual $2$-forms, we
get the corresponding vertical vectors at a spinor $s$ by Clifford
multiplication since the curvature of the connection is
anti-self-dual.
\end{rmk}
\begin{rmk} \label{orientationsrmk}
{\em A note on orientations.} With our chosen convention for the \SP\
$4$-form \Ph, the natural local model for this structure is the
negative spinor bundle over $\mathbb R^4$. With the opposite choice of
orientation, we would be working with the positive spinor
bundle. See~\cite{Korient} for more about sign conventions and
orientations. As we are working only on $S^4$ in this paper, it does
not make a difference.
\end{rmk}
Let $e_1,e_2,e_3,e_4$ be an oriented orthonormal frame for $S^4$
with horizontal lifts to the total space $\spm (S^4)$ denoted by
$\bar{e}_i$ with dual $1$-forms $\bar{e}^i$. Let $\cf^1, \cf^2,
\cf^3, \cf^4$ be the corresponding oriented orthonormal basis for
the fibres. Then (dropping the wedge product symbols for clarity),
the form $\Ph$ can be written as
\begin{eqnarray} \label{BSPhieq}
\Ph & = & u^4 \bar{e}^1 \bar{e}^2 \bar{e}^3 \bar{e}^4 + u^2 v^2
(\bar{e}^1 \bar{e}^2 - \bar{e}^3 \bar{e}^4) (\cf^1 \cf^2 - \cf^3
\cf^4) \\ \nonumber & & {}+ u^2 v^2 (\bar{e}^1 \bar{e}^3 -
\bar{e}^4 \bar{e}^2) (\cf^1 \cf^3 - \cf^4 \cf^2) \\ \nonumber & &
{}+ u^2 v^2 (\bar{e}^1 \bar{e}^4 - \bar{e}^2 \bar{e}^3) (\cf^1
\cf^4 - \cf^2 \cf^3) + v^4 \cf^1 \cf^2 \cf^3 \cf^4
\end{eqnarray}

Now let $\Sigma^2 \subset S^4$ be an oriented surface equipped with
the induced metric and let $(e_1, e_2, \nu_1, \nu_2)$ be an
oriented adapted frame along the surface. That is, $(e_1, e_2)$
are orthonormal tangent vectors and $(\nu_1, \nu_2)$ are
orthonormal normal vectors to the surface. We are interested in
the operator
\begin{equation*}
\Gamma = \gamma(e^1 \wedge e^2) = \pm \gamma(\nu^1 \wedge \nu^2)
\qquad \text{on } \sppm
\end{equation*}
acting on spinors. The operator $\Gamma$ leaves $\sppm$
invariant and it is easily seen that $\Gamma$ is well defined
globally and is independent of the local frame. Moreover $\Gamma$
is a skew-hermitian operator satisfying $\Gamma ^2 = -1$. The
eigenspace decomposition of $\spm$ with respect to $\Gamma$
defines a natural splitting of the spinor bundle $\spm$ restricted
to the surface: $ \rest{\spm}{\Sigma}
\cong \spm^{+} \oplus \spm^{-}$, where
\begin{equation*}
\spm^{\pm} = \{ s \in \spm \, | \,\, \Gamma(\Sigma) = \gamma(e^1
\wedge e^2) s = \pm i \,s \}
\end{equation*}
The two bundles $ \spm^{+}$ and $ \spm^{-}$ are complex line bundles
and are orthogonal to each other. We want to determine necessary and
sufficient conditions on the second fundamental form of $\Sigma$ for
the total space of these bundles to be Cayley submanifolds with
respect to the Bryant-Salamon $\SP$-structure on $M^8$.

\begin{thm} \label{cayleythm}
The total space of either rank $2$ bundle $\spm^{\pm }$ over $\Sigma$
is a Cayley submanifold of $\spm (S^4)$ if and only the immersion
$\Sigma \subset S^4$ is {\em minimal}.
\end{thm}
\begin{proof}
We show every tangent space to the total space of $\spm^{+}$ is a
Cayley subspace of the corresponding tangent space to $\spm
(S^4)$. The proof for $\spm^{-}$ is identical.

Let $ \dot \Gamma $ denote the covariant derivative of the operator
$\Gamma$ along the surface.  Since $\Gamma^2 = -1$, we have: $ \Gamma
\dot \Gamma + \dot \Gamma \Gamma = 0 $, so $\Gamma$ and $\dot \Gamma$
anti-commute and $\dot \Gamma$ interchanges the two eigenspaces of
$\Gamma$. Differentiating the eigenvalue equation $\Gamma s = i s$,
we get $(\Gamma - i) \dot s = - \dot \Gamma s $ and hence
\begin{equation*}
\dot s = - \frac{1}{2} \, i \, \dot \Gamma s
\end{equation*}

Now at a fixed point on $S^4$ let $s_1$ be a unit spinor in the fibre
$\spm^{+}$. Then $s_2 = \Gamma s_1 = i s_1$ is another unit spinor in
$\spm^{+}$ orthogonal to $s_1$. Therefore the fibres of the negative
spinor bundle at a point are given by $t_1 s_1 + t_2 s_2$ where $t_1,
t_2 \in \R$. Thus the following four vectors form a basis of the
tangent space at $t_1 s_1 + t_2 s_2$ of $\spm^{+}$:
\begin{eqnarray*}
E_1 & = & \bar{e}_1 - \frac{i}{2} \,t_1 \, \nabla_{e_1} (\Gamma)
(s_1) - \frac{i}{2} \,t_2 \, \nabla_{e_1} (\Gamma)
(s_2) 
\\ E_2 & = & \bar{e}_2 - \frac{i}{2} \,t_1 \, \nabla_{e_2} (\Gamma)
(s_1) - \frac{i}{2} \,t_2 \, \nabla_{e_2} (\Gamma)
(s_2) \\ F_1 & = & s_1 \\ F_2 &=& s_2 = i s_1
\end{eqnarray*}
where the bar denotes the horizontal lift and the $\nabla_{e_i} (
\Gamma) (s_j)$ are vertical vectors which can be
expressed in terms of the second fundamental form of the
submanifold as we now describe.

Using the adapted frame $(e_1, e_2, \nu_1, \nu_2)$, we have at a
given point (recall we are always using normal coordinates)
\begin{eqnarray*}
\nabla_{e_k} \Gamma & = & \left( \gamma( \nabla_{e_k} e^1)
\gamma(e^2) + \gamma(e^1)\gamma (\nabla_{e_k} e^2) \right) \\ & = &
-A_{k1}^1 \gamma (\nu^1 \wedge e^2) - A_{k1}^2 \gamma (\nu^2 \wedge
e^2) - A_{k2}^1 \gamma (e^1 \wedge \nu^1) - A_{k2}^2 \gamma (e^1
\wedge \nu^2)
\end{eqnarray*}
where we have used the notation $A_{kj}^l = \langle \nabla_{e_k} e_j ,
\nu_l \rangle $. Note that the operators $\gamma(e^j \wedge \nu^l)$
all anti-commute with $\Gamma = \gamma(e^1 \wedge e^2)$ as expected
and hence they permute the two subbundles $\spm^{\pm}$.  Let $\cf^1$
be the $1$-form dual to the vertical tangent vector $\cf_1$ which
corresponds to the spinor $s_1$. Then one can check easily that
$\cf^2, \cf^3, \cf^4$ correspond to the spinors $s_2 =
\frac{\omega_1}{2} \cdot s_1, s_3 = \frac{\omega_2}{2} \cdot s_1, s_4
= \frac{\omega_3}{2} \cdot s_1$, respectively. It can also be checked
that $\gamma(e^1) \gamma(\nu^1) = \gamma(e^2) \gamma(\nu^2)$ and
$\gamma(e^1) \gamma(\nu^2) = - \gamma(e^2) \gamma(\nu^1)$, since we
are on the negative spinor bundle so Clifford multiplication by
$-\gamma(e^1 e^2 \nu^1 \nu^2)$ is equal to $-1$. Using all these facts
the tangent vectors can be expressed as
\begin{eqnarray*}
E_1 & = & \bar{e}_1 + \frac{t_1}{2} \left( (-A^1_{11} - A^2_{12})
\cf_3 + (-A^2_{11} + A^1_{12}) \cf_4 \right) \\ & & {} +
\frac{t_2}{2} \left( (A^2_{11} - A^1_{12}) \cf_3 + (-A^1_{11} -
A^2_{12}) \cf_4 \right) \\ E_2 & = & \bar{e}_2 + \frac{t_1}{2}
\left( (-A^1_{12} - A^2_{22}) \cf_3 + (-A^2_{12} + A^1_{22}) \cf_4
\right) \\ & & {}+ \frac{t_2}{2} \left( (A^2_{12} - A^1_{22}) \cf_3
+ (-A^1_{12} - A^2_{22}) \cf_4 \right)\\ F_1 & = & \cf_1 \\ F_2 &
= & \cf_2 
\end{eqnarray*}
In order to check that the space spanned by $E_1, E_2, F_1, F_2$ is
Cayley, we need to check the vanishing of the $\wedge^2_7$ form
$\eta$ from Proposition~\ref{cayleyprop} using the explicit form of
$\Ph$ in~\eqref{BSPhieq}. Recall that from~\eqref{BS8metriceq} we
have that ${\bar{e}_k}^{\flat} = u^2
\bar{e}^k$ and ${\cf_k}^{\flat} = v^2 \cf^k$. Then (again
omitting the wedge product symbols), one can tediously compute
that
\begin{eqnarray*}
\eta & = & 2 u^2 v^2 \left( t_1(A^1_{11} + A^1_{22}) - t_2(A^2_{11} +
A^2_{22})\right) \left( \bar{e}^1 \cf^3 - \bar{e}^2 \cf^4 -
\bar{e}^3 \cf^1 + \bar{e}^4 \cf^2 \right) \\ & & {}+ 2 u^2 v^2 \left(
t_2(A^1_{11} + A^1_{22}) + t_1(A^2_{11} + A^2_{22})\right) \left(
\bar{e}^1 \cf^4 + \bar{e}^2 \cf^3 - \bar{e}^3 \cf^2 - \bar{e}^4
\cf^1 \right)
\end{eqnarray*}
which clearly vanishes for all $t_1, t_2$ if and only if $\Sigma$
is minimal in $S^4$.
\end{proof}
An obvious example again in this case is to take an equatorial
$S^2$ sitting inside $S^4$, which is totally geodesic. Then there
exist two different real rank $2$ vector bundles over this $S^2$
which are Cayley with respect to the Bryant-Salamon metric on $\spm
(S^4)$. In fact by the results of Bryant~\cite{Br3}, any genus
Riemann surface may be immersed in $S^4$ as a minimal surface, and
hence we can find Cayley submanifolds of $\spm (S^4)$ which are
rank $2$ bundles over any possible compact surface.

\section{Local Intersections of Calibrated Submanifolds}
\label{intersectionsec}

In this section we make some remarks about possible uses of these
constructions to study the local intersections of compact
calibrated submanifolds in a compact manifold with special
holonomy. In~\cite{M} McLean studied the local moduli spaces of
compact calibrated submanifolds. One of his observations was the
following.
\begin{thm}[McLean~\cite{M}] \label{mcleanthm}
Let $X$ be a compact calibrated submanifold of a manifold $M$ with
special holonomy. A small neighbourhood of $X$ in $M$ is naturally
isomorphic to a small neighbourhood of the zero section of the
normal bundle $N(X)$ of $X$ in $M$. We also have the following
explicit identifications of $N(X)$ for the various cases of
calibrations:
\begin{equation*}
\begin{matrix} \text{\underline{Calibration}} &
\text{\underline{Normal Bundle $N(X)$ is isomorphic to}} \\
\text{special Lagrangian} & \text{Cotangent bundle $T^*(X)$
(intrinsic)} \\ \text{coassociative} & \text{Bundle of
anti-self-dual $2$-forms $\wedge^2_- (X)$ (intrinsic)} \\
\text{associative} & \text{twisted spinor bundle
$\spi \otimes_{\Q} E$ over $X$ (non-intrinsic)} \\ \text{Cayley} &
\text{twisted negative spinor bundle $\spm \! \otimes_{\Q} F$ over
$X$ (non-intrinsic)} \\ \end{matrix}
\end{equation*}
where $E$ and $F$ are some explicitly described quaternionic line
bundles.
\end{thm}

Now in all the explicit non-compact manifolds with complete
metrics of special holonomy that we have been discussing in this
paper, the base of the bundle (the zero section), is an example of
a calibrated submanifold. (In fact the zero section is always
rigid with respect to deformations through calibrated submanifolds
by the results of McLean~\cite{M}.) Explicitly, $S^n$ is special
Lagrangian in $T^* (S^n)$ with respect to the Stenzel metric,
$\C \PR^2$ is coassociative in $\wedge^2_- (\C \PR^2)$ with
respect to the Bryant-Salamon metric, and so on. The ambient
manifolds in all cases are complete versions of the local
neighbourhoods described in Theorem~\ref{mcleanthm}. This is
immediate for the special Lagrangian and coassociative cases. In
the case of $S^4$, McLean shows that the quaternionic line bundle
$F$ is trivial in this case so the normal bundle is isomorphic to
$\spm (S^4)$, which is the ambient space of the
complete Bryant-Salamon \SP\ metric. Finally, there is also a
complete \G\ metric on $\spi (S^3)$ that was discovered by Bryant
and Salamon~\cite{BS}. We do not discuss this metric in the
current paper because the calculations are almost identical to the
$\spm (S^4)$ case, but see~\cite{IKM} for some brief remarks on
this metric. The zero section $S^3$ is associative in $\spi(S^3)$,
and the quaternionic line bundle $E$ mentioned in
Theorem~\ref{mcleanthm} is again trivial in this case.

Hence we see that these non-compact manifolds (at least near the
zero section) are good local models for a small neighbourhood of a
rigid, compact calibrated submanifold. Furthermore, one can check
that in all these cases the fibres of the vector bundle total
space are also calibrated submanifolds. The fibres are examples of
calibrated submanifolds which intersect the base calibrated
submanifold in only a point. However, the calibrated submanifolds
which we constructed in Sections~\ref{stenzelsec}
and~\ref{exceptionalsec} were defined as sub-bundles of the total
space restricted to a submanifold of the base. These calibrated
submanifolds interesect the base calibrated submanifold in a
surface in the exceptional cases, and in submanifolds of many
different possible dimensions in the special Lagrangian case.

From the characterizations of calibrated submanifolds in terms of
cross product structures and calibrating forms in
Section~\ref{calibreviewsec}, one can deduce that
(non-singular) calibrated submanifolds can only intersect in
submanifolds of certain allowable dimensions. For instance, since
an associative $3$-plane is closed under the cross product, two
associative $3$-planes can only intersect in $0,1,$ or $3$
dimensions. This is because if they intersect in $2$ dimensions
spanned by orthogonal vectors $e_1$ and $e_2$, the fact that they
are both associative means that must also both contain the third
direction $e_1 \times e_2$. Now because coassociative $4$-planes
are orthogonal complements to associative $3$-planes, one can use
a similar argument to show that two coassociative submanifolds can
only intersect in $0,2,$ or $4$ dimensions. Similarly since Cayley
$4$-planes are closed under the triple cross product $X$, it is
easy to deduce that they too can only intersect in $0,2,$ or $4$
dimensions. Finally, consider the local model of $\R^n \subset
\C^n$ of a special Lagrangian of phase $0$ in $\C^n$, with
coordinated $z^j = x^j + i y^j$. Then the real $n$-plane with
coordinates $(x^1, \ldots, x^p, i y^{p+1}, \ldots, i y^n)$ is a
$\operatorname{U}(n)$ rotation of $\R^n$ with determinant $i^{n-p}$
and hence is special Lagrangian in $\C^n$ with phase $i^{n-p}$,
and intersects $\R^n$ in $p$ dimensions. Thus we have essentially
shown the following.
\begin{prop} \label{intersectionsprop}
Let $X_1$ and $X_2$ be two non-singular calibrated submanifolds of
a manifold $M$ with special holonomy. Suppose that $X_1$ and
$X_2$ intersect at some point $\mathbf x$, and that in a
neighbourhood $U$ of $\mathbf x$ the intersection $X_1 \cap X_2$
is not just the point $x$ and not all of $X_1 \cap U$ (and
equivalently not all of $X_2 \cap U$.) Then we must have:
\begin{equation*}
\begin{matrix} \text{\underline{Calibration}} &
\text{\underline{Intersection of $X_1$ and $X_2$ near $\mathbf x$
must be}} \\ \text{special Lagrangian} & \text{$p$-dimensional,
when phases of $X_1$, $X_2$ differ by $i^{n-p}$} \\
\text{coassociative} & \text{a surface ($2$-dimensional)} \\
\text{associative} & \text{a curve ($1$-dimensional)} \\
\text{Cayley} & \text{a surface ($2$-dimensional)} \\ \end{matrix}
\end{equation*}
\end{prop}
The constructed calibrated submanifolds in this paper all
intersect the base (zero section) calibrated submanifold in
precisely the dimensions expected by
Proposition~\ref{intersectionsprop}. (Compare
Remark~\ref{slagintersectionsrmk}.) Furthermore, our constructions
required strong conditions on the intersection with the base,
thought of as an isometrically immersed submanifold of the base.
Based on this evidence, it is natural to ask the following question:
\begin{quest} \label{intersectionsquest}
Let $X_1$ and $X_2$ be two compact calibrated submanifolds of a
compact manifold $M$ with special holonomy. Recall that both $X_1$ and
$X_2$ inherit induced Riemannian metrics $g_1$ and $g_2$ from $M$,
respectively. Suppose that $X_1$ and $X_2$ intersect at some point
$\mathbf x$, and that in a neighbourhood $U$ of $\mathbf x$ the
intersection $X_1 \cap X_2$ is not just the point $x$ and not all of
$X_1 \cap U$ (and equivalently not all of $X_2 \cap U$.) Then is it
true that we must have:
\begin{itemize}
\item if $X_1$ and $X_2$ are special Lagrangian, with phases
differing by $i^{n-p}$, then the local intersection of $X_1$ and
$X_2$ near $\mathbf x$ is a $p$-dimensional submanifold, which is
an austere immersion with respect to $(X_1,g_1)$ or $(X_2,g_2)$. 
\item if $X_1$ and $X_2$ are coassociative, then the local
intersection of $X_1$ and $X_2$ near $\mathbf x$ is a
$2$-dimensional surface, which is a properly oriented isotropic
(that is, negative superminimal) immersion with respect to
$(X_1,g_1)$ or $(X_2,g_2)$.
\item if $X_1$ and $X_2$ are associative, then the local
intersection of $X_1$ and $X_2$ near $\mathbf x$ is a
$1$-dimensional curve, which is a geodesic (minimal) immersion with
respect to $(X_1,g_1)$ or $(X_2,g_2)$.
\item if $X_1$ and $X_2$ are Cayley, then the local
intersection of $X_1$ and $X_2$ near $\mathbf x$ is a
$2$-dimensional surface, which is a minimal immersion with respect
to $(X_1,g_1)$ or $(X_2,g_2)$.  
\end{itemize}
\end{quest}
We are currently investigating this question. A related problem is the
following. In symplectic geometry, a neighbourhood of a Lagrangian
submanifold $X$ in a symplectic manifold $M$ is naturally identified
with a neigbourhood of the zero section in $T^*(X)$. It would be
useful to have similar neighbourhood theorems in the case of
calibrated submanifolds, describing the Ricci-flat metric on the
ambient space to a certain order of approximation. Topologically, this
was done by McLean~\cite{M}.

It would also be useful to discover to what extent these bundle
constructions of calibrated submanifolds generalize to other
explicitly known metrics. There is a wealth of new explicit examples
of \G\ and \SP\ metrics, for example, that have been recently
discovered by physicists.  (See~\cite{CGLP1, CGLP2}, and the
references therein.)

\section{Conclusion} \label{conclusionsec}

Besides the possible applications to the study of intersections of
calibrated submanifolds discussed in Section~\ref{intersectionsec},
there are several other future directions to explore. It would be
interesting to study the possible singularities that can occur in such
examples. It should be noted that even when the submanifold over which
we build our calibrated sub-bundle is only immersed in the base, with
self-intersections, the resulting calibrated submanifold which we
construct is in fact embedded. It is also worth studying how these
calibrated submanifolds can be deformed. This would require extending
the work of McLean~\cite{M} to the case of non-compact calibrated
submanifolds. Some study has been made of deformations of non-compact
asymptotically conical~\cite{J4,L3,Ma,P} or asymptotically
cylindrical~\cite{JS} calibrated submanifolds. This of course is
closely related to the possible non-existence of other kinds of
calibrated submanifolds built as bundles over the same submanifold,
discussed at the end of Section~\ref{intersectionsec}. It may be that
the only way to deform our constructed calibrated submanifolds through
calibrated submanifolds would be to deform the base of the
sub-bundle. For example, the moduli space of associative $3$-folds
near a fixed associative submanifold $L$ which is a rank $1$ line
bundle over a minimal surface $\Sigma$ as constructed in
Section~\ref{exceptionalsec} may be just those which arise via the
same construction by deforming the minimal surface inside the base,
through minimal surfaces. These moduli of course always exist as
possible deformations, the only question being whether or not there
are any others.


\begin{thebibliography}{99}

\bibitem{An}H. Anciaux, {\em Special Lagrangian Submanifolds in
the Complex Sphere}, preprint, arXiv:math.DG/0311288.

\bibitem{AW}M. Atiyah and E. Witten, {\em M-Theory Dynamics on a
Manifold of \G\ Holonomy}, Adv. Theor. Math. Phys. {\bf 6} (2003),
1-106.

\bibitem{Br3}R.L. Bryant, {\em Conformal and Minimal Immersions of
Compact Surfaces into the Four-Sphere}, J. Diff. Geom. {\bf 17}
(1982), 455-473.

\bibitem{Br1}R.L. Bryant, {\em Submanifolds and Special Structures
on the Octonions}, J. Diff. Geom. {\bf 17} (1982),
185-232.

\bibitem{Br2}R.L. Bryant, {\em Metrics with Exceptional Holonomy},
Ann. of Math. {\bf 126} (1987), 525-576.

\bibitem{Br4}R.L. Bryant, {\em Some remarks on the geometry of
austere manifolds}, Bol. Soc. Brasil. Mat. (N.S.) {\bf 21} (1991),
133-157.

\bibitem{BS}R.L. Bryant and S.M. Salamon, {\em On the Construction
of Some Complete Metrics with Exceptional Holonomy}, Duke Math. J.
{\bf 58} (1989), 829-850.

\bibitem{Ca}E. Calabi, {\em M\'etriques k\"ahl\'eriennes et
fibr\'es holomorphes}, Ann. Sci. Ecole Norm. Sup. (4) {\bf 12}
(1979), 269-294.

\bibitem{CGLP1}M. Cveti\v c, G.W. Gibbons, H. L\"u, and C.N.
Pope, {\em Ricci-flat Metrics, Harmonic Forms, and Brane
Resolutions}, Comm. Math. Phys. {\bf 232} (2003), 457-500.

\bibitem{CGLP2}M. Cveti\v c, G.W. Gibbons, H. L\"u, and C.N.
Pope, {\em Hyper-K\"ahler Calabi Metrics, $L^2$ Harmonic Forms,
Resolved M2-branes, and $AdS_4$/$CFT_3$ Correspondence}, Nucl.
Phys. B {\bf 617} (2001), 151-197.

\bibitem{DF}M. Dajczer and L.A. Florit, {\em A Class of Austere
Submanifolds}, Illinois J. Math. {\bf 45} (2001), 735-755.

\bibitem{Dst}A. Dancer and I.A.B. Strachan, {\em Einstein Metrics
on Tangent Bundles of Spheres}, Classical Quantum Gravity {\bf 19}
(2002), 4663-4670.

\bibitem{DS}A. Dancer and A. Swann, {\em Hyperk\"ahler Metrics of
Cohomogeneity One}, J. Geom. Phys. {\bf 21} (1997), 218-230.

\bibitem{ES}J. Eells and S. Salamon, {\em Twistorial Construction
of Harmonic Maps of Surfaces into Four-Manifolds}, Ann. Scuola
Norm. Sup. Pisa. Cl. Scil. (4) {\bf 12} (1985), 589-640. 

\bibitem{GPP}G.W. Gibbons, D.N. Page, and C.N. Pope, {\em Einstein
Metrics on $S^3$, $\mathbb R^3$, and $\mathbb R^4$ bundles},
Commun. Math. Phys. {\bf 127} (1990), 529-553.

\bibitem{GS}S. Gukov and J. Sparks, {\em M-Theory on Spin(7)
Manifolds}, Nucl. Phys. B {\bf 625} (2002), 3-69.

\bibitem{H}R. Harvey, {\em Spinors and Calibrations}, Academic
Press, San Diego, 1990.

\bibitem{HL2}R. Harvey and H.B. Lawson, {\em A Constellation of
Minimal Varieties Defined Over the Group $G_2$},
Lecture Notes in App. Math. {\bf 48} (1979), 167-187.

\bibitem{HL}R. Harvey and H.B. Lawson, {\em Calibrated
Geometries}, Acta Math. {\bf 148} (1982), 47-157. 

\bibitem{Hi}N. Hitchin, {\em The Geometry of Three-Forms in Six
and Dimensions}, J. Diff. Geom. {\bf 55} (2000), 547-576.

\bibitem{IKM}M. Ionel, S. Karigiannis, and M. Min-Oo, {\em Bundle
Constructions of Calibrated Submanifolds in $\R^7$ and $\R^8$},
arXiv: math.DG/0408005, submitted for publication.

\bibitem{J1}D.D. Joyce, {\em Compact Manifolds with Special
Holonomy}, Oxford University Press, Oxford, 2000.

\bibitem{J2}D.D. Joyce, {\em Lectures on Special Lagrangian
Geometry}, arXiv: math.DG/0111111.

\bibitem{J3}D.D. Joyce, {\em The Exceptional Holonomy Groups and
Calibrated Geometry}, arXiv: math.DG/0406011.

\bibitem{J4}D.D. Joyce, {\em Special Lagrangian Submanifolds with
Isolated Conical Singularities. II. Moduli Spaces}, Ann. Global
Anal. Geom. {\bf 25} (2004), 303-352.

\bibitem{JS}D.D. Joyce and S. Salur, {\em Deformations of
Asymptotically Cylindrical Coassociative Submanifolds with Fixed
Boundary}, arXiv: math.DG/0408137.

\bibitem{K}S. Karigiannis, {\em Deformations of \G\ and \SP\
structures}, Canad. J. Math., to appear.

\bibitem{Korient}S. Karigiannis, {\em A note on signs and orientations
in \G\ and \SP\ geometry}; posted on the arXiv.

\bibitem{LL}J.H. Lee and N.C. Leung, {\em Instantons and Branes in
Manifolds with Vector Cross Product}, arXiv: math.DG/0402044.

\bibitem{LW}N.C. Leung and X.W. Wang, {\em Intersection Theory of
Coassociative submanifolds in \G\ manifolds and Seiberg-Witten
invariants}, arXiv: math.DG/0401419.

\bibitem{L1}J. Lotay, {\em Constructing Associative $3$-folds by
Evolution Equations}, arXiv: math.DG/0401123.

\bibitem{L2}J. Lotay, {\em $2$-Ruled Calibrated $4$-folds in
$\mathbb R^7$ and $\mathbb R^8$}, arXiv: math.DG/0401125.

\bibitem{L3}J. Lotay, {\em Deformation Theory of Asymptotically
Conical Coassociative $4$-folds}, arXiv: math.DG/0411116.

\bibitem{Ma}S.P. Marshall, {\em Deformations of Special Lagrangian
Submanifolds}, Oxford DPhil. Thesis, 2002.

\bibitem{M}R.C. McLean, {\em Deformations of Calibrated
Submanifolds}, Comm. Anal. Geom. {\bf 6} (1998), 705-747.

\bibitem{P}T. Pacini, {\em Deformations of Asymptotically Conical
Special Lagrangian Submanifolds}, Pacific J. Math., to appear.

\bibitem{Sa2}S.M. Salamon, {\em Harmonic and Holomorphic Maps},
Lecture Notes in Math. {\bf 1164} (1984), 161-224.

\bibitem{Sa}S.M. Salamon, {\em Riemannian Geometry and Holonomy
Groups}, Longman Group UK Limited, Harlow, 1989.

\bibitem{St}M.B. Stenzel, {\em Ricci-flat Metrics on the
Complexification of a Compact Rank One Symmetric Space},
Manuscripta Math. {\bf 80} (1993), 151-163.

\bibitem{SYZ}A. Strominger, S.T. Yau, and E. Zaslow, {\em Mirror
Symmetry is T-Duality}, Nuclear Physics B {\bf 479} (1996),
243-259.

\bibitem{Sz}R. Sz\"oke, {\em Complex Structures on
Tangent Bundles of Riemannian Manifolds}, Math. Ann. {\bf 291}
(1991), 409-428.

\bibitem{W}S.H. Wang, {\em On the Lifts of Minimal Lagrangian
Submanifolds}, arXiv: math.DG/0109214.

\end{thebibliography}
\end{document}